\documentclass[12pt]{article}
\usepackage[latin1]{inputenc}
\usepackage{amssymb,amsmath,amsfonts}
\usepackage{xcolor}
\newtheorem{theorem}{Theorem}
\newtheorem{proposition}[theorem]{Proposition}
\newtheorem{lemma}[theorem]{Lemma}

\newtheorem{remark}[theorem]{Remark}

\newcommand{\R}{\mathbb{R}}
\newcommand{\U}{\mathcal{U}}
\newcommand{\C}{\mathbb{C}}
\newcommand{\spa}{\mbox{span\,}}
\newcommand{\rank}{\mbox{rank }}
\newcommand{\po}{{\hspace*{-1ex}}{\bf .  }}

\newcommand{\nap}{\nabla^{\perp}}
\newcommand{\nab}{\tilde\nabla}
\def\lp{{\langle\!\langle}}\vspace{2ex}
\def\rp{{\rangle\!\rangle}}
\def\<{{\langle}}
\def\>{{\rangle}}
\def\Sal{{\cal S}}
\def\J{{\cal J}}

\def\a{\alpha}
\def\be{\begin{equation} }
\def\ee{\end{equation} }

\def\nap{\nabla^\perp}
\def\proof{\noindent\emph{Proof: }}
\def\qed{\ifhmode\unskip\nobreak\fi\ifmmode\ifinner
\else\hskip5 pt \fi\fi\hbox{\hskip5 pt \vrule width4 pt
height6 pt  depth1.5 pt \hskip 1pt }}

\newcommand\blfootnote[1]{
\begingroup
\renewcommand\thefootnote{}\footnote{#1}%
\addtocounter{footnote}{-1}
\endgroup}
\begin{document}

\title{Real Kaehler submanifolds in codimension\\ up to four}
\author{S. Chion and M. Dajczer}
\date{}
\maketitle

\begin{abstract}  
Let $f\colon M^{2n}\to\mathbb{R}^{2n+4}$ be an isometric immersion 
of a Kaehler manifold of complex dimension $n\geq 5$ into 
Euclidean space with complex rank at least $5$ everywhere. 
Our main result is that, along each connected component of 
an open dense subset of $M^{2n}$, either $f$ is holomorphic 
in $\mathbb{R}^{2n+4}\cong\mathbb{C}^{n+2}$ or it is in a unique way a 
composition $f=F\circ h$ of isometric immersions. In the latter 
case, we have that $h\colon M^{2n}\to N^{2n+2}$ is holomorphic 
and $F\colon N^{2n+2}\to\mathbb{R}^{2n+4}$ belongs to the class, 
by now quite well understood, of non-holomorphic Kaehler submanifold 
in codimension two. Moreover, the submanifold $F$ is minimal if 
and only if $f$ is minimal. 
\end{abstract}
\blfootnote{\textup{2020} \textit{Mathematics Subject Classification}:
53B25, 53B35.}

By a \emph{real Kaehler submanifold} $f\colon M^{2n}\to\R^{2n+p}$ 
we mean an isometric immersion of a connected Kaehler manifold 
$(M^{2n},J)$ of complex dimension $n\geq 2$ into Euclidean space 
with codimension $p$. Throughout this paper, it is assumed that $f$ 
is locally substantial, which means that the image of $f$ restricted 
to any open subset of  $M^{2n}$ does not lie inside a proper affine 
subspace of $\R^{2n+p}$. Moreover, if the codimension is even then $f$, 
when restricted to any open subset of $M^{2n}$, is not allowed 
to be holomorphic in $\R^{2n+2q}\cong\C^{n+q}$. Notice that conditions 
that yield that $f$ is holomorphic have been given in \cite{CCD} and 
that $f$ is just minimal in \cite{CD}.
\vspace{1ex}

The study of real Kaehler submanifolds has acquired increasing 
relevance since the pioneering work by Dajczer and Gromoll in \cite{DG1}. 
Clearly, the main motivation for their study has been that when minimal 
these submanifolds  enjoy many of the  feature properties 
of minimal surfaces. For instance, if simply-connected they admit 
an associated one-parameter family of non congruent isometric minimal 
submanifolds with the same Gauss map and can be realized as the 
real part of its holomorphic representative.  Moreover, they are 
pluriharmonic (sometimes called pluriminimal) submanifolds and,
in particular, they are austere submanifolds. Furthermore, as seen in 
the Appendix, there are several cases when a classification 
is reached through a Weierstrass type representation. For a 
partial account of results and references on this subject of 
research we refer to \cite{DT}.\vspace{1ex}

It is well-known that the second fundamental form 
$\a\colon TM\times TM\to N_f M$ of a real Kaehler submanifold 
$f\colon M^{2n}\to\R^{2n+p}$ with codimension $p=1$ or $p=2$ has 
necessarily a large kernel.  The complex dimension of that kernel
is measured by way of the notion of (complex) rank.
\vspace{1ex}

The \emph{rank} $\varrho_f(x)$ of $f$ at $x\in M^{2n}$ is given by  
$2\varrho_f(x)=2n-\dim\Delta(x)\cap J\Delta(x)$ where 
\begin{equation*}
\Delta(x)=\{Y\in T_xM\colon \a(X,Y)=0\;\;\mbox{for any}\;\;X\in T_xM\}
\end{equation*}
is known as the relative nullity subspace of $f$ at $x\in M^{2n}$.
\vspace{1ex}

As recalled below, outside a flat point the rank is $\varrho_f=1$ if 
the codimension is $p=1$, and that it is $\varrho_f\leq 2$ at any point 
if $p=2$. The situation is quite different for submanifolds in higher 
codimension, in part due to the presence of compositions of isometric 
immersion. For instance, already for $p=3$ we have that if 
$F\colon N^{2n+2}\to\R^{2n+3}$ is a real Kaehler hypersurface and 
$g\colon M^{2n}\to N^{2n+2}$ an holomorphic submanifold, then the 
composition of isometric immersions $f=F\circ g\colon M^{2n}\to\R^{2n+3}$ 
may have any rank since there is no bound under holomorphicity. Thus, 
in the search of local classifications of real Kaehler submanifolds 
in higher codimension than two, but still low, a necessary step is 
to provide conditions that impose the existence of a composition.

At this time there is substantial knowledge about the local real 
Kaehler submanifolds that are free of flat points and lie in 
codimension of at most four. On one hand, the ones in codimension 
one or two are quite well understood. On the other hand, for 
the higher codimensions three and four it turns out that 
under a proper rank assumption the submanifold has to be a 
composition as the one for $p=3$ discussed above.  
\vspace{1ex}

The non flat real Kaehler hypersurfaces $f\colon M^{2n}\to\R^{2n+1}$, 
$n\geq 2$, have been locally classified by Dajczer and Gromoll 
\cite{DG1} by way of the so called Gauss parametrization in terms 
of a pseudoholomorphic surface in the $2n$-dimensional round sphere 
and any smooth function on the surface. This was made possible 
because in this case the rank is $\varrho_f=1$; see Theorem $15.14$ in 
\cite{DT} for a more detailed proof of this classification. It turns 
out that the hypersurface is minimal if and only if the function on 
the surface is an eigenvector of the Laplacian for the eigenvalue 
$2$. For this special case, there is a Weierstrass type parametrization 
given by Hennes \cite{He}.

For real Kaehler submanifolds $f\colon M^{2n}\to\R^{2n+2}$, $n\geq 3$,
Dajczer \cite{Da} showed that $\dim\Delta(x)\geq 2n-4$ at any 
$x\in M^{2n}$. For a discussion of the classification for the 
real Kaehler submanifolds that lie in codimension two in terms of  
its rank see the Appendix in this paper.
\vspace{1ex}

We say that real Kaehler submanifold $f\colon M^{2n}\to\R^{2n+p}$ 
admits a \emph{Kaehler extension} if there exist a (maybe flat) real 
Kaehler submanifold $F\colon N^{2n+2\ell}\to\R^{2n+p}$, $\ell\geq 1$, 
and a holomorphic isometric embedding $j\colon M^{2n}\to N^{2n+2\ell}$ 
such that $f=F\circ j$. 
\vspace{1ex}

For $f\colon M^{2n}\to\R^{2n+3}$, $n\geq 4$, it was stated by Dajczer 
and Gromoll \cite{DG4} that if $\dim\Delta<2n-6$ everywhere,
then there exists an open dense subset of $M^{2n}$ such that, along 
each connected component, the submanifold admits a unique Kaehler 
extension to a real Kaehler hypersurface. Unfortunately, the existence 
part of the proof in \cite{DG4}  depends on an algebraic lemma that went 
unproven. But that result is correct for the specific non symmetric bilinear 
form given by ($1$) in that paper, as follows from Lemma~\ref{unique} 
or from Lemma $1$ in \cite{YZ}. A nice observation due to Yan and Zheng 
\cite{YZ} is that this result should hold under a weaker assumption 
on the rank of the submanifold. In this paper we provide a proof of this.

The case of real Kaehler submanifolds $f\colon M^{2n}\to\R^{2n+4}$,
$n\geq 5$, is treated by Yan and Zheng in \cite{YZ} and \cite{YZ2}. 
It is proved in \cite{YZ} that, if the rank of $f$ 
satisfies $\varrho_f>4$ everywhere, there exists an open dense 
subset of $M^{2n}$ such that $f$ restricted to each connected component 
admits a Kaehler extension $F\colon N^{2n+2}\to\R^{2n+4}$. 
In the present paper, by way of an alternative approach, we are able 
to complement in some directions the Main Theorem in \cite{YZ}. 
In particular, we show that their result is correct in spite of a 
rather minor inaccuracy in the proof. In fact, their statement that 
the vanishing of a shape operator on one normal direction amounts 
to a reduction of the codimension of the submanifold does not hold. 
It is just elementarily false that in this case there is a
reduction of the  codimension and, of course, there is no statement 
in Spivak that proves such a claim.  
Nevertheless, it turns out that even in this case their theorem is 
correct since, in this situation, the submanifold is holomorphic 
inside a flat submanifold as the one given in the non minimal 
case $(i)-(a)$ in the Appendix. 

The main achievement of this paper is to prove that the Kaehler 
extension is \emph{unique}, up to reparametrizations, in opposition 
to a assertion in \cite{YZ}. Consequently, the submanifold $f$ 
is minimal if and only if its extension $F$ is also minimal.
What made our proof of uniqueness possible is that the way
we construct the extensions is somehow more restricted than 
in \cite{YZ}; see Remark \ref{remark}.

\begin{theorem}\po\label{main}
Let $f\colon M^{2n}\to\R^{2n+p}$, $3\leq p\leq 4$ and $n>p$, be 
a real Kaehler submanifold whose rank satisfies $\varrho_f>p$ everywhere.
Then there exists an open dense subset of $M^{2n}$ such that the
restriction of $f$ to each connected components admits a unique
Kaehler extension $F\colon N^{2n+2}\to\R^{2n+p}$. 
Moreover, the  rank of $F$ is constant $\varrho_F\leq p-2$ and $F$ is 
a minimal submanifold if and only if $f$ is minimal.
\end{theorem}

Under the rank assumption required in the above result, the manifold
$M^{2n}$ is free of points where all the sectional curvatures vanish.
In fact, by a classical result going back to Cartan, at such a point we 
have $\dim\Delta\geq 2n-p$ and hence $\varrho_f\leq p$.
\vspace{1ex}

Theorem \ref{main} is sharp even if  the submanifold is asked to be 
isometrically complete. For instance, the extrinsic product immersion 
of two complete minimal ruled submanifolds lying in codimension two 
classified in \cite{DG3} has rank four.
\vspace{1ex}

From the Appendix of this paper it follows that the geometric 
options for the extension $F$, and hence for $f$, in Theorem \ref{main} 
not to be minimal are quite limited.
\vspace{1ex}

We observe that Yan and Zheng in \cite{YZ} made
a very interesting and rather challenging conjecture: If the 
codimension is $p\leq 11$ then Kaehler extensions always exist 
for real Kaehler submanifold $f\colon M^{2n}\to\R^{2n+p}$, $n>p$, 
if the rank satisfies $\varrho_f>p$ everywhere. In this respect, 
the results in the first section of this paper, that hold for any 
codimension, should be useful. Finally, in \cite{CD2} it is 
shown why in the Yang and Zheng conjecture limiting the codimension 
to $11$ is essential.  This was done proving that the structure 
of the second fundamental form until that codimension is the one
expected and this by an argument that fails beyond that codimension. 
As for the full conjecture, we intend to give an answer in a 
forthcoming paper that makes use of some results from the present 
one.

\section{A class of Kaehler extensions} 

The goal of this section is to establish a set of conditions 
for a real Kaehler submanifold to admit a Kaehler extension 
of a certain type. The result achieved holds regardless of 
the size the codimension and should be of use for further 
applications.
\vspace{1ex}

We first introduce some notations and definitions.  
Let $\gamma\colon V\times V\to W$ be a bilinear form between 
real vector spaces. The image of $\gamma$ is the vector 
subspace of $W$ given by
\begin{equation*}
\mathcal{S}(\gamma)
=\spa\{\gamma(X,Y)\;\mbox{for all}\; X,Y\in V\}
\end{equation*}
whereas the (right) nullity of $\gamma$ is the vector subspace 
of $V$ defined by
\begin{equation*}
\mathcal{N}(\gamma)=\{Y\in V\colon\gamma(X,Y)=0
\;\mbox{for all}\;X\in V\}.
\end{equation*}
Henceforth $f\colon M^{2n}\to\R^{2n+p}$ stands for a 
real Kaehler submanifold and
\begin{equation*}
\a\colon TM\times TM\to N_f M
\end{equation*}
its second fundamental form. If $P\subset N_fM(x)$ is 
a vector subspace and if $\a_P$ is the $P$-component of the second 
fundamental form $\a$ at $x\in M^{2n}$ then the \emph{complex kernel} 
of $\a_P$ is the tangent vector subspace
\begin{equation*}
\mathcal{N}_c(\a_P)
=\mathcal{N}(\a_P)\cap J\mathcal{N}(\a_P).
\end{equation*}

Throughout this section $L$ denotes a normal vector subbundle of 
(real) rank $2\ell>0$ that satisfies $L(x)\subset\Sal(\a(x))$ 
everywhere and is endowed with the induced metric and vector bundle 
connection. Moreover, it is assumed that $L$ carries an isometric 
complex structure $\J\in\Gamma(Aut(L))$, that is, a vector bundle 
isometry that satisfies $\J^2=-I$. Furthermore, it is required the 
complex tangent vector subspaces $D(x)=\mathcal{N}_c(\a_{L^\perp}(x))$ 
to have constant (even) dimension and thus form a holomorphic tangent 
subbundle denoted by $D$.
\vspace{1ex}

In the sequel the pair $(L,\J)$ is required to satisfy  the 
following conditions:
\begin{itemize}
\item[($\mathcal{C}_1$)] The complex structure 
$\J\in\Gamma(Aut(L))$ is parallel,
that is, 
\begin{equation*}
(\nap_X\J\eta)_L=\J(\nap_X\eta)_L \;\;\mbox{for any}\;\;X\in\mathfrak{X}(M)
\;\mbox{and}\;\eta\in\Gamma(L)
\end{equation*}
and  the second fundamental form of $f$ satisfies
\be\label{Jsecondf}
\J\a_L(X,Y)=\a_L(X,JY)\;\;\mbox{for any}\;\; X,Y\in\mathfrak{X}(M)
\ee
or, equivalently, the shape operators verify 
$A_{\J\eta}= J\circ A_\eta=-A_\eta\circ J$ for $\eta\in\Gamma(L)$.
\item[($\mathcal{C}_2$)] The subbundle $L$ is parallel along $D$ in 
the normal connection of $f$, that is,
\begin{equation*}
\nap_Y\eta\in\Gamma(L)\;\;\mbox{for any}\;\;Y\in\Gamma(D)\;\mbox{and}\;  
\eta\in\Gamma(L).
\end{equation*}
\end{itemize}

We observe that the  subbundle $L$ is
necessarily proper in $N_1=\Sal(\a)$ along
any open subset $U\subset M^{2n}$ where
the latter has constant rank since, otherwise, we would have from the
condition $(\mathcal{C}_2)$ that $N_1=N_fM$ and then $f$ would be 
holomorphic along $U$ as established by Proposition \ref{holo} below. 
\vspace{1ex}

Let the vector bundle $TM\oplus L$ over $M^{2n}$ be endowed 
with the complex structure $\hat\J\in\Gamma(Aut(TM\oplus L))$ 
defined by 
\be\label{hatJ}
\hat\J(X+\eta)=JX+\J\eta.
\ee 
Condition $(\mathcal{C}_1)$ easily gives that $\hat\J$ is 
parallel in the induced vector bundle connection defined by
$\hat\nabla_X(Y+\eta)=(\nab_X(Y+\eta))_{TM\oplus L}$
where $\nab$ denotes the Euclidean connection in the ambient 
space $\R^{2n+p}$. That is, we have
\be\label{parallelhatj}
(\nab_X\hat\J(Y+\eta))_{TM\oplus L}
=\hat\J((\nab_X(Y+\eta))_{TM\oplus L})
\ee 
for any $X,Y\in\mathfrak{X}(M)$ and $\eta\in\Gamma(L)$. 

\begin{proposition}\po\label{charD} 
The distribution $D=\mathcal{N}_c(\a_{L^\perp})$ is integrable.
\end{proposition}

\proof The Codazzi equation $(\nabla_X^\perp\a)(S,T)
=(\nabla_S^\perp\a)(X,T)$ and the condition $(\mathcal{C}_2)$ give
\be\label{codhojapluri}
(\nabla_X^\perp\a(S,T))_{L^\perp}+\a_{L^\perp}(X,\nabla_ST)=0
\ee
for any $X\in\mathfrak{X}(M)$ and $S,T\in\Gamma(D)$. 

From \eqref{codhojapluri} it follows that  
$[S,T]\in\mathcal{N}(\a_{L^\perp})$ for any $S,T\in\Gamma(D)$.
On the other hand, from \eqref{Jsecondf} we have
$\a(S,JT)=\a(JS,T)$ for any $S,T\in\Gamma(D)$.
This and \eqref{codhojapluri} give
\begin{equation*}
\a_{L^\perp}(X,J[S,T])=\a_{L^\perp}(X,\nabla_SJT-\nabla_TJS)
=(\nabla_X^\perp(\a(T,JS)-\a(S,JT)))_{L^\perp}=0
\end{equation*}
for any $X\in\mathfrak{X}(M)$ and $S,T\in\Gamma(D)$.
Thus also  $J[S,T]\in\mathcal{N}(\a_{L^\perp})$ and
hence $[S,T]\in D$ for any $S,T\in\Gamma(D)$. 
\vspace{2ex}\qed

Proposition \ref{charD} gives that $M^{2n}$ carries an holomorphic 
foliation. Let $i\colon\Sigma\to M^{2n}$ denote the inclusion of 
the leaf $\Sigma$ through $x\in M^{2n}$ and let 
$g\colon\Sigma\to\R^{2n+p}$ be the isometric immersion $g=f\circ i$.

We have that
$$
\nab_ST=f_*(\nabla_ST)_D+f_*(\nabla_ST)_{D^\perp}+\a^f(S,T)
$$
for any $S,T\in\Gamma(D)$. Hence
\be\label{secondfundg}
\a^g(y)(S,T)=f_*(i(y))(\nabla_{i_*S}i_*T)_{D^\perp}
+\a^f(i(y))(i_*S,i_*T)
\ee
for any $S,T\in\mathfrak{X}(\Sigma)$.  

We have that $g=f\circ i\colon\Sigma\to\R^{2n+p}$ 
satisfies $g(\Sigma)\subset f_*T_xM\oplus L(x)$. To see this,  
observe that the normal bundle of $g$ splits  orthogonally as 
$N_g\Sigma=i^*(f_*D^\perp\oplus L\oplus L^\perp)$ and that the 
condition $(\mathcal{C}_2$) gives that the vector subbundle 
$i^*L^\perp$ of $N_g\Sigma$ is constant in $\R^{2n+p}$. Thus, 
we may also see $g(\Sigma)$ as a submanifold of
$\R^{2n+2\ell}=f_*T_xM\oplus L(x)\subset\R^{2n+p}$
when it is convenient.

\begin{proposition}\po\label{charD2} The submanifolds 
$g\colon\Sigma\to\R^{2n+2\ell}$ are holomorphic. Moreover, for 
any given $g$ the map $\psi\colon\Sal(\a^g)\to\Sal(\a^f|_{D\times D})$ 
defined by
$$
\psi(\a^g(S,T))=\a^f(i_*S,i_*T)
$$
is an isomorphism. 
\end{proposition}

\proof Let $\bar\J$ be the complex structure on $\R^{2n+2\ell}$ 
induced by $\hat\J$. Then 
$$
\bar\J g_*T=\hat\J f_*i_*T=f_*Ji_*T=f_*i_*J|_{T\Sigma}T
=g_*J|_{T\Sigma}T
$$
for any $T\in\mathfrak{X}(\Sigma)$, and hence $g$ is holomorphic.

From \eqref{secondfundg}  the map $\psi$ is surjective. 
To prove injectivity, 
let $\delta=\sum_{j=1}^k\a^g(S_j,T_j)$ for
$S_j,T_j\in\mathfrak{X}(\Sigma)$ satisfy
$\psi(\delta)=\sum_{j=1}^k\a^f(i_*S_j,i_*T_j)$=0.
From  \eqref{Jsecondf} we have
$$
\sum_{j=1}^k\a^f(i_*S_j,i_*JT_j)=\J\sum_{j=1}^k\a^f(i_*S_j,i_*T_j)=0.
$$
We obtain from \eqref{codhojapluri} that
$$
\sum_{j=1}^k\a_{L^\perp}(X,\nabla_{S_j}i_*T_j)
=-\sum_{j=1}^k(\nabla_X^\perp\a^f(i_*S_j,i_*T_j))_{L^\perp}=0,
$$
and, similarly, that $\sum_{j=1}^k\a_{L^\perp}(X,J\nabla_{S_j}i_*T_j)=0$.
Hence $\sum_{j=1}^k\nabla_{S_j}i_*T_j\in\Gamma(D)$ and we conclude from 
\eqref{secondfundg} that $\delta=0$.\vspace{2ex}\qed

In the sequel the pair $(L,\J)$ is assumed to  satisfy  the 
additional condition:
\vspace{2ex} 

\noindent $(\mathcal{C}_3)$ $L=\Sal(\a|_{D\times D})$
at any point of $M^{2n}$.

\begin{lemma}\po\label{existence}
Let $\pi\colon\Lambda\to M^{2n}$ be the $\hat\J$-invariant 
vector subbundle of $TM\oplus L$ defined~by 
\be\label{Lambda}
\Lambda=\spa\{(\nabla_ST)_{D^\perp}+\a^f(S,T)\colon S,T\in\Gamma(D)\}.
\ee
Then $\rank\Lambda=2\ell$ and $\Lambda\cap TM=0$. Moreover, we have that
\be\label{derivative}
\nab_X\lambda\in f_*TM\oplus L
\ee 
for any $\lambda\in\Gamma(\Lambda)$ and $X\in\mathfrak{X}(M)$.
\end{lemma}

\proof From Proposition \ref{charD2} it follows that 
$\rank\Sal(\a^g)=\rank i^*L$. If $\Sigma$ is the leaf of $D$ that 
contains $x\in M^{2n}$ we have from \eqref{secondfundg} that 
$\Lambda(i(x))=\Sal(\a^g)(x)$.  Thus $\Lambda$ has constant
dimension at each $x\in M^{2n}$ and hence is a subbundle.   
If $\lambda\in\Lambda\cap TM$ 
it follows from Proposition \ref{charD2} and \eqref{secondfundg} 
that $\lambda=0$ and hence $\Lambda\cap TM=0$.

We obtain from \eqref{codhojapluri} that
$$
(\nab_X(f_*(\nabla_ST)_{D^\perp}+\a^f(S,T)))_{L^\perp}=0
$$
for any $S,T\in\Gamma(D)$ as we wished.\vspace{2ex}\qed

\begin{lemma}\po A real Kaehler submanifold 
$f\colon M^{2n}\to\R^{2n+p}$ is minimal if and only 
if it is pluriharmonic, that is, if 
\be\label{pluriharmonic}
\a(JX,Y)=\a(X,JY)\;\;\mbox{for any}\;\; X,Y\in \mathfrak{X}(M)
\ee
or, equivalently, the shape operators verify 
$J\circ A_\xi=-A_\xi\circ J$ for any $\xi\in N_fM$.
\end{lemma}

\proof This is Theorem $1.2$ in \cite{DR} or Theorem $15.7$ 
in \cite{DT}.\qed

\begin{theorem}\po\label{develop} Let $f\colon M^{2n}\to\R^{2n+p}$ 
be an embedding and let $\pi\colon\Lambda^{2\ell}\to M^{2n}$ be
the vector subbundle of  $TM\oplus L^{2\ell}$ defined by \eqref{Lambda}. 
Let $N^{2n+2\ell}$ be an open neighborhood of the $0$-section  
$j\colon M^{2n}\to N^{2n+2\ell}$ of $\Lambda^{2\ell}$ such that 
the map $F\colon N^{2n+2\ell}\to\R^{2n+p}$ defined by 
$$
F(\lambda)=f(\pi(\lambda))+\lambda
$$
is an embedding. Then $F$ is a Kaehler extension of $f=F\circ j$ 
and its second fundamental form satisfies 
$\mathcal{N}_c(\a^F)=D\oplus\Lambda$ at any point of $M^{2n}$ 
and thus, in particular, its rank is $2\varrho_F=2n-\rank D$. 
Moreover, the submanifold $F$ is minimal if and only if $f$ is minimal.
\end{theorem}

\proof If $\lambda_0\in N^{2n+2\ell}$ and $x_0=\pi(\lambda_0)$
let $\xi\in\Gamma(\Lambda)$  be such that $\xi(x_0)=\lambda_0$.
Then
$$
\nab_XF(\lambda_0)=f_*(x_0)X+\nab_X\xi(x_0)\;\;
\mbox{if}\;\;X\in T_{x_0}M.
$$
Thus $T_{\lambda}N=T_{\pi(\lambda)}M\oplus L(\pi(\lambda))$
and $N_FN(\lambda)=L^\perp(\pi(\lambda))$ for any 
$\lambda\in N^{2n+2\ell}$.

We first show that $F$ is a real Kaehler submanifold. 
Let $J^N\in\Gamma(Aut(TN))$ be the complex structure defined 
by $J^N(\lambda)=\hat{\J}(\pi(\lambda))$ where $\hat{\J}$ 
is given by \eqref{hatJ}. Thus $J^N$ is constant along the 
fibers of $\Lambda$.
Then $J^N$ is parallel with respect to  the Levi-Civita 
connection $\hat{\nabla}$ on $N^{2n+2\ell}$ since by 
\eqref{parallelhatj} we have that
$$
\hat{\nabla}_XJ^N(Y+\xi)=
(\tilde{\nabla}_X\hat{\J}(Y+\xi))_{TM\oplus L}
=\hat{\J}(\tilde{\nabla}_X(Y+\xi))_{TM\oplus L}
=J^N\hat{\nabla}_X (Y+\xi)
$$
for any $X,Y\in\mathfrak{X}(M)$ and $\xi\in\Gamma(L)$.

We show next that $\mathcal{N}_c(\a^F)=D\oplus\Lambda$.
The condition $(\mathcal{C}_2)$ gives
$$
\nab_S\eta=-A_\eta^f S+\nabla_S^\perp\eta
=\nabla_S^\perp\eta\in\Gamma(L^\perp)
$$
for any $S\in\Gamma(D)$ and $\eta\in\Gamma(L^\perp)$.  
We obtain $D\oplus\Lambda\subset\mathcal{N}(\a^F)$. Since 
$D$ is $J$-invariant and $\Lambda$ is $\hat\J$-invariant
then $D\oplus\Lambda$ is $J^N$-invariant. Thus
$D\oplus\Lambda\subset\mathcal{N}_c(\a^F)$.
For the other inclusion, we have to verify that if
$Z\in\mathfrak{X}(M)\cap\mathcal{N}_c(\a^F)$ then
$Z\in\Gamma(D)$. Since
$$
{}^F\nabla_Z^\perp\eta=\nab_Z\eta
=-A_\eta^fZ+{}^f\nabla_Z^\perp\eta \;\;\mbox{and}\;\;
{}^F\nabla_{JZ}^\perp\eta=\nab_{JZ}\eta
=-A_\eta^fJZ+{}^f\nabla_{JZ}^\perp\eta
$$
for any $\eta\in\Gamma(L^\perp)$ then
$A_\eta^fZ=A_\eta^fJZ=0$.

Assume that $f$ is minimal. We first prove the following fact:
\be\label{developmin}
\J(\nabla_X^\perp\eta)_L+(\nabla_{JX}^\perp\eta)_L=0\;\;
\mbox{for any}\;X\in\mathfrak{X}(M)\;\mbox{and}\;\eta\in\Gamma(L^\perp).
\ee
Let $\lambda\in\Gamma(\Lambda)$ be such that $(\lambda)_L
=\J(\nabla_X^\perp\eta)_L+(\nabla_{JX}^\perp\eta)_L$.
Using first that \eqref{pluriharmonic} holds and then 
\eqref{derivative} at the end of the argument, we obtain
\begin{align*}
\|(\lambda)_L\|^2&=
\<\lambda,\J(\nabla_X^\perp\eta)_L+(\nabla_{JX}^\perp\eta)_L\>
=\<\lambda,-JA_\eta X+\J(\nabla_X^\perp\eta)_L-A_\eta JX
+(\nabla_{JX}^\perp\eta)_L\>\\
&=\<\lambda,\hat\J(\nab_X\eta)_{TM\oplus L}
+(\nab_{JX}\eta)_{TM\oplus L}\>
=\<\nab_X\hat\J\lambda,\eta\>-
\<\nab_{JX}\lambda,\eta\>=0.
\end{align*}
Since $\nab_X\eta=-A_\eta^fX+\nabla_X^\perp\eta$ then
$A_\eta^FX=A_\eta^fX-(\nabla_X^\perp\eta)_L$.
Using \eqref{developmin} we have
\begin{align*}
A_\eta^FJ^NX
&=A_\eta^fJX-(\nabla_{JX}^\perp\eta)_L
=-JA_\eta^fX+\J(\nabla_X^\perp\eta)_L
=-\hat\J(A_\eta^fX-(\nabla_X^\perp\eta)_L)\\
&=-J^NA^F_\eta X
\end{align*}
for any $\eta\in\Gamma(L^\perp)$ and $X\in\mathfrak{X}(M)$.
Since we have seen that $\mathcal{N}_c(\a^F)=D\oplus\Lambda$ 
then $\a^F(\delta,\xi)=\a^F(\delta,J^N\xi)=0$
holds for any $\delta\in TN$ and $\xi\in\Lambda$, and
hence $F$ is a minimal immersion.

Assume that $F$ is minimal. Since $J$ is holomorphic then
$\a^j(JS,T)=\a^j(S,JT)$ for any $S,T\in\mathfrak{X}(\Sigma)$. 
Now since $f=F\circ j$ then
\begin{align*}
\a^f(S,JT)
&=F_*\a^j(S,JT)+\a^F(j_*S,j_*JT)=F_*\a^j(JS,T)+\a^F(j_*S,J^Nj_*T)\\
&=F_*\a^j(JS,T)+\a^F(J^Nj_*S,j_*T)=F_*\a^j(JS,T)+\a^F(j_*JS,j_*T)\\
&=\a^f(JS,T)
\end{align*}
and hence $f$ is a minimal submanifold.\qed

\begin{remark}\po\label{remark}
{\em In Proposition 2 in \cite{YZ} the extension is obtained 
by means of a developable ruling $L$ whereas here Theorem
\ref{develop} uses the subbundle $\Lambda$ given by \eqref{Lambda}.
This is more restricted since $\Lambda$
is a special case of a canonical developable ruling as defined
in \cite{YZ}. In fact, our $\Lambda$ is $J$-invariant and 
that may not be the case of $L$. 
}\end{remark}

\section{The proof of Theorem \ref{main}}

The proof of Theorem \ref{main} requires several results. The 
first one holds for any codimension and is of independent interest.

\begin{proposition}\po\label{holo}
Let $f\colon M^{2n}\to\R^{2n+p}$ be an isometric immersion of a 
Kaehler manifold. Assume that $N_1(x)=\Sal(\a(x))$ satisfies 
$N_1(x)=N_fM(x)$ at any $x\in M^{2n}$ and there is an
isometry $\J\in\Gamma(Aut(N_fM))$ such that 
\be\label{thecond}
\J\a(X,Y)=\a(X,JY)\;\;\mbox{for any}\;\; X,Y\in\mathfrak{X}(M).
\ee
Then $p$ is even and $f$ is holomorphic.
\end{proposition}

\proof From \eqref{thecond} we have that the $\J$ is 
a complex structure which we claim to be parallel in the normal 
connection. If we apply $\J$ to the  Codazzi equation 
$(\nabla_X^\perp\a)(Y,Z)=(\nabla_Y^\perp\a)(X,Z)$  subtract 
$(\nabla_X^\perp\a)(Y,JZ)=(\nabla_Y^\perp\a)(X,JZ)$ and then
use \eqref{thecond} we obtain
\be\label{kappaeq}
\mathcal{K}(X)\a(Y,Z)=\mathcal{K}(Y)\a(X,Z)
\ee
where $\mathcal{K}(X)\in\Gamma(\text{End}(N_fM))$ is the 
skew-symmetric tensor defined by
$$
\mathcal{K}(X)\eta
=\J\nabla_X^\perp\eta-\nabla_X^\perp\J\eta.
$$
If we  denote
$$
\<\mathcal{K}(X)\a(Y,Z),\a(S,T)\>=(X,Y,Z,S,T)
$$
then by \eqref{kappaeq} 
and since $\mathcal{K}(X)$ is skew-symmetric, we obtain
\begin{align*}
(X,Y,Z,S,T)
&=-(X,S,T,Y,Z)=-(S,X,T,Y,Z)=(S,Y,Z,X,T)=(Y,S,Z,X,T)\\
&=-(Y,X,T,S,Z)=-(T,X,Y,S,Z)=(T,S,Z,X,Y)=(Z,S,T,X,Y)\\
&=-(Z,X,Y,S,T)=-(X,Y,Z,S,T)=0
\end{align*}
for any $X,Y,Z,S,T\in\mathfrak{X}(M)$. Because $N_1=N_fM$
everywhere then $\mathcal{K}(X)=0$ for any $X\in\mathfrak{X}(M)$, 
and the claim has been proved. 

If follows from the claim that $\J$ is constant along the submanifold 
and hence $J\oplus\J$ extends to a complex structure on $\R^{2n+p}$, 
still denoted by $\J$, such that $\J\circ f_*=f_*\circ\J$, and 
therefore $f$ is holomorphic. \vspace{2ex}\qed

Let $f\colon M^{2n}\to\R^{2n+p}$ be a real Kaehler submanifold 
and let $N_fM(x)\oplus N_fM(x)$ be endowed with the inner product
of signature $(p,p)$ given by
$$
\lp(\xi_1,\xi_2),(\eta_1,\eta_2)\rp
=\<\xi_1,\eta_1\>-\<\xi_2,\eta_2\>.
$$
We call a bilinear form  
$\varphi\colon T_xM\times T_xM\to N_fM(x)\oplus N_fM(x)$ 
\emph{flat} if 
$$
\lp\varphi(X,Y),\varphi(Z,T)\rp-\lp\varphi(X,T),\varphi(Z,Y)\rp=0
$$
for any $X,Y,Z,T\in T_xM$.  

\begin{lemma}\po\label{flatforms} The bilinear form 
$\gamma\colon T_xM\times T_xM\to N_1(x)\oplus N_1(x)$ 
defined by
\be\label{gamma}
\gamma(X,Y)=(\a(X,Y),\a(X,JY))
\ee 
is flat.
\end{lemma}

\proof It is well-known that the curvature tensor of a Kaehler 
manifold $M^{2n}$ satisfies $R(X,Y)=R(JX,JY)$ and $R(X,Y)JZ=JR(X,Y)Z$ 
for any $X,Y,Z\in T_xM$. Then a roughly short straightforward computation 
making use of this as well as the Gauss equation of $f$ gives the flatness.
\vspace{2ex}\qed 

A vector subspace $V\subset N_fM(x)\oplus N_fM(x)$ is called 
a degenerate space if $V\cap V^\perp\neq 0$ and if otherwise 
nondegenerate. 

\begin{lemma}\po\label{prop5} 
Let $f\colon M^{2n}\to\R^{2n+p}$, $p<n$, be an isometric immersion 
of a Kaehler manifold and let $P\subset N_1(x)$ be a vector subspace 
of dimension $\dim P\leq 5$. Assume that the bilinear form 
$\gamma_P\colon T_xM\times T_xM\to P\oplus P$ defined by
$$
\gamma_P(X,Y)=(\a_P(X,Y),\a_P(X,JY))
$$
is flat and the vector space $\Sal(\gamma_P)$
is nondegenerate. Then 
$\dim\mathcal{N}_c(\a_P(x))\geq 2n-2\dim P$.  

In particular, if $p\leq 5$ and the flat bilinear form $\gamma$ in 
\eqref{gamma} satisfies that $\Sal(\gamma)$ is a nondegenerate
vector space then $\varrho_f(x)\leq p$. 
\end{lemma}

\proof The proof follows from Proposition $5$ in \cite{CCD2}.\qed

\begin{lemma}\po\label{prop5alt}
Let $f\colon M^{2n}\to\R^{2n+p}$, $p<n$ and $2\leq p\leq 5$ be an 
isometric immersion of a Kaehler manifold with rank $\varrho_f>p$ 
everywhere. \vspace{1ex}

\noindent $(i)$ At any $x\in M^{2n}$ there is a subspace 
$L(x)\subset N_1(x)$ of
$\dim L(x)=2\ell(x)>0$, an isometry $\J(x)\in Aut(L(x))$ such that
\be\label{sffpuntual}
\J(x)\a_{L(x)}(X,Y)=\a_{L(x)}(X,JY)\;\;\mbox{for any}\;\; X,Y\in T_xM,
\ee
the subspace $\Sal(\gamma_{L^\perp(x)})$ is nondegenerate and
$\dim\mathcal{N}_c(\a_{L(x)^\perp})\geq 2n-2p+4\ell(x)$.
\vspace{1ex}

\noindent $(ii)$  On each connected component of an open dense 
subset $U\subset M^{2n}$ we have  $\ell(x)=\ell$ and
$\dim\mathcal{N}_c(\a_{L(x)^\perp})$ are constant, the subspaces 
$L(x)$ form a normal vector subbundle $L\subset N_1$ of 
$\rank L=2\ell$, there is an isometry $\J\in\Gamma(Aut(L))$
that satisfies \eqref{Jsecondf} and there is a tangent subbundle 
$\mathcal{N}_c(\a_{L^\perp})$ such that
$\rank\mathcal{N}_c(\a_{L^\perp})\geq 2n-2p+4\ell$.
\end{lemma}

\proof By Lemma \ref{prop5} the bilinear form $\gamma$
given by \eqref{gamma} satisfies that $\Sal(\gamma)$ is a
degenerate subspace at any $x\in M^{2n}$, that is, we have that
$\U(x)=\Sal(\gamma)\cap(\Sal(\gamma))^\perp\neq 0$.  If
$(\xi,\eta)\in\U(x)$ then also $(\eta,-\xi)\in\U(x)$
since
$$
\lp\gamma(X,Y),(\eta,-\xi)\rp=\lp\gamma(X,JY),(\xi,\eta)\rp
$$
for any $X,Y\in T_xM$. Then $\dim\U(x)=2\ell(x)>0$.

We have that $\pi_1(\U(x))=\pi_2(\U(x))$ where
$\pi_j\colon N_1(x)\oplus N_1(x)\to N_1(x)$, $j=1,2$,
is the projection onto the $j$-th component. Since
$\pi_j|_{\U(x)}$, $j=1,2$, is injective then $L(x)=\pi_j(\U(x))$
satisfies $\dim L(x)= 2\ell(x)$. We have 
$\U(X)\subset\Sal(\gamma_{L(x)})\subset L(x)\oplus L(x)$ and 
thus, by dimension reasons, we obtain $\U(X)=\Sal(\gamma_{L(x)})$.  
Then
\be\label{flatness}
\<\a_{L(x)}(X,Y),\a_{L(x)}(Z,T)\>=\<\a_{L(x)}(X,JY),\a_{L(x)}(Z,JT)\>
\ee
for any $X,Y,Z,T\in T_xM$. Thus, there is an isometry
$\J(x)\in Aut(L(x))$ such that \eqref{sffpuntual} holds.
We have that $\gamma$ is flat and  \eqref{flatness} just says 
that also
$$
\gamma_{L(x)}(X,Y)=(\a_{L(x)}(X,Y),\a_{L(x)}(X,JY)).
$$
is flat. Since $\gamma(x)=\gamma_{L(x)}+\gamma_{L^\perp(x)}$
then also $\gamma_{L^\perp(x)}$ is flat. Having that subspace
$\Sal(\gamma_{L^\perp(x)})$ is nondegenerate then Lemma \ref{prop5} 
gives $\dim\mathcal{N}_c(\a_{L^\perp(x)})\geq 2n-2p+4\ell(x)$.
\vspace{1ex}

Let $U\subset M^{2n}$ be an open sense subset  such that $\ell(x)$
and $\dim\mathcal{N}_c(\a_{L^\perp(x)})$ are constant on each 
connected component. Along any component the subspaces $\U(x)$
form a vector bundle and thus also the subspaces $L(x)$ do since
they possess equal dimension.  Finally, that $\J\in\Gamma(Aut(L))$
is smooth follows from \eqref{sffpuntual}. \vspace{2ex}\qed

For codimension $p=2$ the following result generalizes the one
in \cite{Da}. A proof should also follow from the arguments in \cite{YZ}.

\begin{theorem}\po\label{cod2}
Let $f\colon M^{2n}\to\R^{2n+2}$, $n\geq 3$ be an isometric immersion 
of a Kaehler manifold. If the rank is $\varrho_f>2$ everywhere then 
$f$ is an holomorphic submanifold.  
\end{theorem}

\proof Lemma \ref{prop5alt} gives that $N_1=N_fM$  
and an isometry $\J\in\Gamma(Aut(N_fM))$ satisfying
$$
\J\a(X,Y)=\a(X,JY)
$$
for any $X,Y\in\mathfrak{X}(M)$. Then Proposition \ref{holo} 
yields that $f$ is holomorphic. \qed

\begin{lemma}\po\label{unique}
Let $f\colon M^{2n}\to\R^{2n+p}$, $3\leq p\leq 4$ and $n>p$, be 
a real Kaehler submanifold with rank $\varrho_f>p$ everywhere. 
Then along each connected component, say $U_0$, of an open dense subset
of $M^{2n}$ there exists a pair $(L,\J)$ such that 
$L\subset N_1|_{U_0}$ is a vector subbundle of rank $2$ and 
$\J\in\Gamma(\text{Aut}(L))$ an isometric complex structure that 
satisfies~\eqref{Jsecondf}.

Moreover, the subspaces $D(x)=\mathcal{N}_c(\a_{L^\perp}(x))$ form a 
tangent vector subbundle such that $\rank D\geq 2n-2p+4$. Furthermore, 
if $(L',\J')$ is a pair along $U_0$ such that $L'\subset N_1$ is
a vector subbundle and the isometry $\J'\in\Gamma(\text{Aut}(L'))$ 
satisfies \eqref{Jsecondf} then $(L,\J)=(L',\J')$.
\end{lemma}

\proof  By Lemma \ref{prop5alt} on each connected component of an open
dense subset $U$ of $M^{2n}$ there is a normal vector subbundle
$L\subset N_1$ and an isometry $\J\in\Gamma(Aut(L))$ such that 
$\rank L=2$ if $p=3$ and either $\rank L=2$ or $\rank L=4$ if $p=4$.
Moreover, $D=\mathcal{N}_c(\a_{L^\perp})$ satisfies 
$\dim D\geq 2n-2p+2\,\rank L$.
If $\rank L=4$ then $p=4$ and $L=N_1=N_fM$. Hence \eqref{Jsecondf}
holds for $L=N_fM$ and Proposition \ref{holo} yields that $f$ is 
holomorphic, which has been excluded. Therefore, $\rank L=2$ on each 
connected component of $U$ and $\rank D\geq 2n-2p+4$.  

We prove the uniqueness part of the statement. From Proposition 
\ref{holo} we obtain $\rank L'=2$. By assumption, we have  
that $\gamma_{L'}(X,Y)=(\a_{L'}(X,Y),\a_{L'}(X,JY))$ satisfies
$$
\lp\gamma_{L'}(X,Y),\gamma_{L'}(Z,T)\rp=0
$$
for any $X,Y,Z,T\in\mathfrak{X}(M)$ and since
$\gamma=\gamma_{L'}+\gamma_{L'^\perp}$ thus $\gamma_{L'^\perp}$
is flat.  We claim that $D'=\mathcal{N}_c(\a_{L'^\perp})$
satisfies $\dim D'\geq 2n-2p+4$. If $\Sal(\gamma_{L'^\perp})$ is
nondegenerate the claim follows from Lemma \ref{prop5}. Thus, 
it suffices to show that 
$\U'=\Sal(\gamma_{L'^\perp})\cap(\Sal(\gamma_{L'^\perp}))^\perp\neq 0$ 
leads to a contradiction. If $(\xi,\bar\xi)\in\U'$ then also 
$(\bar\xi,-\xi)\in\U'$ since
$$
\lp\gamma_{L'^\perp}(X,Y),(\bar\xi,-\xi)\rp=
\lp\gamma_{L'^\perp}(X,JY),(\xi,\bar\xi)\rp=0.
$$
Hence $\dim\U'=2$ and $p=4$. Since
$\U'\subset\Sal(\gamma_{L'^\perp})
\subset(L'^\perp\cap N_1)\oplus(L'^\perp\cap N_1)$
then
\be\label{charU'}
\U'=\Sal(\gamma_{L'^\perp})
\ee
and thus
$\lp\gamma_{L'^\perp}(X,Y),\gamma_{L'^\perp}(S,JT)\rp=0$.
Then there is  a complex structure of the form 
$(\J'\oplus\bar\J)\in\Gamma(Aut(N_fM))$  such that 
$(\J'\oplus\bar\J)\a(X,Y)=\a(X,JY)$.
Being $\pi_j|_{\U'}$ injective then
$\pi_j(\U')=N_1\cap L'^\perp=L'^\perp$ and hence
$L'^\perp\subset N_1$ by \eqref{charU'}. Therefore $N_1=N_fM$ 
and thus $f$ is holomorphic by Proposition \ref{holo}, which 
is not allowed and proves the claim.

Suppose that we have $L\neq L'$. If in case $p=3$ then 
$\dim L^\perp\oplus L'^\perp=2$. Since $\dim D,\dim D'\geq 2n-2$ 
then $\dim\mathcal{N}_c(\a_{L^\perp\oplus L'^\perp})
\geq\dim D\cap D'\geq 2n-4$.  
We have that $L\cap(L^\perp\oplus L'^\perp)\neq 0$ and thus 
\eqref{Jsecondf} gives that $\varrho_f\leq 2$, a contradiction. 
Hence  $p=4$.

If $L^\perp+L'^\perp=N_fM$ we have that 
$D\cap D'\subset\mathcal{N}_c(\a)$ which is not possible since 
it yields $\varrho_f\leq 4$. Hence $\dim(L^\perp+L'^\perp)=3$ and then 
$L\cap(L^\perp+L'^\perp)\neq 0$. Since  \eqref{Jsecondf} is
equivalent to $A_{\J\eta}= J\circ A_\eta$ for any $\eta\in\Gamma(L)$ 
it follows that $D\cap D'\subset\mathcal{N}_c(\a_L)$. 
Hence if $S\in\Gamma(D\cap D')$ then
$$
\a(X,S)=\a_L(X,S)+\a_{L^\perp}(X,S)=0
$$
where the first term on the right-hand-side vanishes
since $S\in\Gamma(D\cap D')$ and the second since $S\in\Gamma(D)$. 
Thus we obtain again that $D\cap D'\subset\mathcal{N}_c(\a)$.\qed

\begin{lemma}\po\label{c1c2}
Let $f\colon M^{2n}\to\R^{2n+p}$, $3\leq p\leq 4$ and $n>p$, be 
a real Kaehler submanifold with rank $\varrho_f>p$ everywhere.
Then any pair $(L,\J)$ given by Lemma \ref{unique} satisfies 
the conditions ($\mathcal{C}_1$), ($\mathcal{C}_2$) and ($\mathcal{C}_3$).
\end{lemma}

\proof That the condition $(\mathcal{C}_1)$ is satisfied 
is trivial.
\vspace{1ex}

For the proof of the condition $(\mathcal{C}_2)$ we first
show that the distribution $D$ is integrable. Given 
$\delta\in\Gamma(L)$ from the Codazzi equation 
$(\nabla_Z A)(\xi;Y)=(\nabla_Y A)(\xi;Z)$ for $\xi=\J\delta$ 
and since $A_{\J\delta}=J\circ A_\delta$ we have
$$
J(\nabla_Z A_\delta Y-A_\delta\nabla_ZY)
-A_{\nabla_Z^\perp\J\delta}Y
=J(\nabla_Y A_\delta Z-A_\delta\nabla_YZ)
-A_{\nabla_Y^\perp\J\delta}Z
$$
for any $Y,Z\in\mathfrak{X}(M)$. From the above Codazzi equation 
for $\xi=\delta$ and the condition $(\mathcal{C}_1)$ we obtain
$$
JA_{(\nabla_Z^\perp\delta)_{L^\perp}}Y
-A_{(\nabla_Z^\perp\J\delta)_{L^\perp}}Y=
JA_{(\nabla_Y^\perp\delta)_{L^\perp}}Z
-A_{(\nabla_Y^\perp\J\delta)_{L^\perp}}Z
$$
for any $Y,Z\in\mathfrak{X}(M)$.  Hence 
\be\label{codl}
\<\a(Y,X),(\nabla_Z^\perp\delta)_{L^\perp}\>
-\<\a(Y,JX),(\nabla_Z^\perp\J\delta)_{L^\perp}\>
=0
\ee 
for any $\delta\in\Gamma(L)$, $Z\in\Gamma(D)$ and 
$X,Y\in\mathfrak{X}(M)$.  

On one hand, the Codazzi equation
$(\nap_{Z_1}\a)(Z,Z_2)=(\nap_{Z_2}\a)(Z,Z_1)$ 
gives
\be\label{liebracket}
\a_{L^\perp}(Z,[Z_1,Z_2])
=(\nap_{Z_1}\a(Z,Z_2)-\nap_{Z_2}\a(Z,Z_1))_{L^\perp}
\ee
for any $Z_1,Z_2\in\Gamma(D)$ and $Z\in\mathfrak{X}(M)$.
On the other hand, the Codazzi equation 
$(\nap_{Z}\a)(Z_i,JZ_j)=(\nap_{Z_i}\a)(Z,JZ_j)$
yields
$$
(\nap_Z\a(Z_i,JZ_j))_{L^\perp}=
(\nap_{Z_i}\a(Z,JZ_j))_{L^\perp}-\a_{L^\perp}(Z,\nabla_{Z_i}JZ_j)
$$
for any $Z_i,Z_j\in\Gamma(D)$ and $Z\in\mathfrak{X}(M)$.
Since $\a(Z_1,JZ_2)=\a(JZ_1,Z_2)$ by \eqref{Jsecondf} then
\be\label{liebracket2}
\a_{L^\perp}(Z,J[Z_1,Z_2])=(\nap_{Z_1}\a(Z,JZ_2)
-\nap_{Z_2}\a(Z,JZ_1))_{L^\perp} 
\ee 
for any $Z_1,Z_2\in\Gamma(D)$ and $Z\in\mathfrak{X}(M)$.

Using first \eqref{liebracket} and \eqref{liebracket2}  and then
\eqref{codl} for $\delta=\a(Z,Z_i)$ we obtain 
\begin{align*}
\<\a_{L^\perp}(Y,X)&,\a_{L^\perp}(Z,[Z_1,Z_2])\>-
\<\a_{L^\perp}(Y,JX),\a_{L^\perp}(Z,J[Z_1,Z_2])\>\\
&=\,\<\a_{L^\perp}(Y,X),\nap_{Z_1}\a(Z,Z_2)\>
-\<\a_{L^\perp}(Y,X),\nap_{Z_2}\a(Z,Z_1)\>\\
&\;\;-\<\a_{L^\perp}(Y,JX),\nap_{Z_1}\a(Z,JZ_2)\>
+\<\a_{L^\perp}(Y,JX),\nap_{Z_2}\a(Z,JZ_1)\>\\
&=0
\end{align*}
for any $Z_1, Z_2\in\Gamma(D)$ and $X,Y,Z\in\mathfrak{X}(M)$. 
Thus we have shown that
$$
\lp\gamma_{L^\perp}(Y,X),\gamma_{L^\perp}(Z,[Z_1,Z_2])\rp=0
$$
for any $Z_1, Z_2\in\Gamma(D)$ and $X,Y,Z\in\mathfrak{X}(M)$.
Since the subspace $\Sal(\gamma_{L^\perp})$ is nondegenerate 
by part $(i)$ of Lemma \ref{prop5alt}, we have 
$\gamma_{L^\perp}(Z,[Z_1,Z_2])=0$ for any $Z_1, Z_2\in\Gamma(D)$
and $Z\in\mathfrak{X}(M)$.  Hence $[Z_1,Z_2]\in\Gamma(D)$ as
we wished.

Since $D$ is integrable then the Codazzi equation
$(\nabla_SA)(\eta, T)=(\nabla_TA)(\eta,S)$ yields
\be\label{codazzi}
A_{(\nabla_S^\perp\eta)_L}T=A_{(\nabla_T^\perp\eta)_L}S
\ee
for any $S,T\in\Gamma(D)$ and $\eta\in\Gamma(L^\perp)$.
If the condition $(\mathcal{C}_2)$ does not hold there is
$S\in\Gamma(D)$ and $\eta\in\Gamma(L^\perp)$ such that
$\mu=(\nabla_S^\perp\eta)_L\neq 0$. If follows from \eqref{codazzi} 
and $\dim D\geq 2n-2p+4$ that $\dim\ker A_\mu\geq 2n-2p+2$
and hence we have by \eqref{Jsecondf} that $\varrho_f\leq p-1$,
which is a contradiction. 

From \eqref{Jsecondf} we have that $\Sal(\a^f|_{D\times D})$ is of 
even dimension. Then $L$ satisfies condition ($\mathcal{C}_3$) since, 
otherwise, we have $\a|_{D\times D}=0$ and hence 
$A_\xi D\subset D^\perp$ if $\xi\in\Gamma(L)$. But then 
$\dim\ker A_\xi|_D\geq 2n-4p+8$ and \eqref{Jsecondf} gives that 
$\varrho_f\leq 2p-4$, which is a contradiction.
\vspace{2ex}\qed

Finally, we are in the condition to prove our main result.
\vspace{2ex}

\noindent\emph{Proof of Theorem \ref{main}:} 
By Lemma \ref{unique} and Lemma \ref{c1c2} there is an open dense 
subset of $M^{2n}$ such that along any connected component, 
say $U$, the submanifold $f|_U$ is an embedding and there
is a unique pair $(L,\J)$,  where $L\subset N_1|_U$ has 
$\rank L=2$,  and a $J$-invariant vector subbundle
$D=\mathcal{N}_c(\a_{L^\perp})$ with $\rank D\geq 2n-2p+4$ such that
conditions $(\mathcal{C}_1)$ to $(\mathcal{C}_3)$ are satisfied.
Then it follows from Theorem \ref{develop} that $f|_U$ admits a Kaehler 
extension $F$ as in the statement.

We now argue for the uniqueness of the Kaehler extension. Let
$F\colon N^{2n+2}\to\R^{2n+p}$, $n\geq p+1$ and $3\leq p\leq 4$, be 
a real Kaehler submanifold such that the tangent vector subspaces
$\Delta^c(z)=\mathcal{N}_c(\a^F(z))$ satisfy that $\dim\Delta^c(z)$ 
is constant and thus form a tangent vector subbundle. In fact,
it is easy to verify that the distribution $\Delta^c$ is integrable 
and that its leaves are totally geodesic submanifolds in $M^{2n}$ 
as well as in $\R^{2n+p}$.

From Lemma \ref{prop5} if $p=3$ and Theorem \ref{cod2} if $p=4$ we have 
$\rank\Delta^c\geq 2n-2p+6$.
Then let $j\colon M^{2n}\to N^{2n+2}$ be an holomorphic submanifold
of $N^{2n+2}$ such that the real Kaehler submanifold
$f=F\circ j\colon M^{2n}\to\R^{2n+p}$ is substantial and satisfies
that $\varrho_f(x)>p$ at any $x\in M^{2n}$.  To conclude the proof, 
we have to show that $F$ is the unique Kaehler extension of $f$
up to a reparametrization.

From Lemmas \ref{unique} and \ref{c1c2} along each connected component,
say $U$, of an open dense subset of $M^{2n}$ there are a unique pair
$(L,\J)$ where $L\subset N_1|_U$ is a vector subbundle
of rank two and the isometric complex structure $\J\in\Gamma(Aut(L))$ 
satisfies
$$
\J\a_L(X,Y)=\a_L(X,JY)\;\;\mbox{for any}\;\; X,Y\in\mathfrak{X}(U).
$$
Moreover, there is a tangent vector subbundle 
$D=\mathcal{N}_c(\a^f_{L^\perp})$ with $\rank D\geq 2n-2p+4$
and the conditions $(\mathcal{C}_1)$ to $(\mathcal{C}_3)$ hold.

Since $f=F\circ j$ then $N_fM=F_*N_jM\oplus j^*N_FN$ and hence
\be\label{sff}
\a^f(X,Y)=F_*\a^j(X,Y)+\a^F(j_*X,j_*Y)
\ee
for any $X,Y\in\mathfrak{X}(M)$. The vector subspace $E\subset T_xM$ 
given by $j_*E=\Delta^c(j(x))\cap j_*T_xM$ satisfies $\dim E\geq 2n-2p+4$
since the codimension of $j_*T_xM$ in $T_{j(x)}N$ is two. 
Being $j$ holomorphic then  $\bar J\a^j(X,Y)=\a^j(X,JY)$ 
for any $X,Y\in\mathfrak{X}(M)$ where $\bar J$ 
is the complex structure of $N^{2n+2}$.  Since $\Delta^c(j(x))\cap j_*T_xM$
is $\bar J$-invariant then $E$ is $J$-invariant and therefore either
$\a^j|_{E\times E}=0$ or $\Sal(\a^j|_{E\times E})=N_jM$.

Suppose that $\a^j|_{E\times E}=0$. If $0\neq\xi\in\Gamma(N_jM)$
and $A^j_\xi$ is the shape operator of $j$ then
$\ker A^j_\xi|_E\subset E\cap\mathcal{N}_c(\a^f)$. In fact, if 
$S\in\ker A^j_\xi|_E$ then $\rank N_jM=2$ gives that 
$S,JS\in\ker A^j_\xi|_E\cap\ker A^j_{\bar J\xi}|_E
\subset\mathcal{N}_c(\a^j)\cap E$. Since $j_*E\subset\Delta^c$ we have 
\mbox{$\a^F(j_*X,j_*S)=0$} if $S\in\Gamma(E)$ and $X\in\mathfrak{X}(M)$. 
Now \eqref{sff} yields $\ker A^j_\xi|_E\subset\mathcal{N}_c(\a^f)$ and 
thus $\dim\mathcal{N}_c(\a^f)\geq 2n-4p+8$. Then $\varrho_f\leq 2p-4$ 
which is a contradiction. It follows that $\Sal(\a^j|_{E\times E})=N_jM$
and we obtain from  \eqref{sff} that
$\Sal(\a^f|_{E\times E})=F_*N_jM\subset N_1$.

The vector subbundle $F_*N_jM\subset N_1$ is endowed with the 
complex structure $F_*\circ\bar{J}|_{N_jM}$ and
$F_*\bar{J}|_{N_jM}\a^j(X,Y)=F_*\a^j(X,JY)$ holds.
Now the uniqueness part of Lemma \ref{unique} yields that
$L=F_*N_jM$ and $\J=F_*\circ\bar{J}|_{F_*N_jM}$.

We claim that $E=D$. Since from \eqref{sff} we have
$\a^f_{L^\perp}(X,Y)=\a^F(j_*X,j_*Y)$ then
$j_*D=\mathcal{N}_c(\a^F|_{j_*TM\times j_*TM})$ and hence
$\Delta^c\cap j_*TM\subset j_*D$.  Given $S\in\Gamma(D)$ then the
condition $(\mathcal{C}_2)$ yields $\nabla_S^\perp F_*\xi\in\Gamma(L)$
for any $\xi\in\Gamma(N_jM)$. Thus $\a^F(j_*S,\xi)=0$ for any
$\xi\in\Gamma(N_jM)$. Since $\a^F(j_*S,j_*X)=\a^f_{L^\perp}(S,X)=0$
for any $X\in\mathfrak{X}(M)$ hence $j_*S\in\Delta^c$ giving the
other inclusion and hence the claim. 

Let $\pi\colon\Lambda\to U$ be the vector subbundle of 
$f_*TU\oplus L$ given by \eqref{Lambda}. Then Lemma~\ref{existence} 
and Theorem \ref{develop} yield the Kaehler extension
$\bar{F}\colon\bar{N}^{2n+2}\subset\Lambda\to\R^{2n+p}$ 
of $f$ defined by $\bar{F}(\lambda)=f(\pi(\lambda))+\lambda$.
To conclude the proof, it remains to show that  $\bar{N}^{2n+2}$ 
can be chosen sufficiently small so that $\bar{F}(\lambda)\in F(N)$ 
if $\lambda=f_*(\nabla_ST)_{D^\perp}+\a^f(S,T)\in\bar{N}^{2n+2}$.
On one hand, from \eqref{sff} and since $j_*D\subset\Delta^c$ 
from the claim we obtain
\begin{align*}
\lambda
&=F_*j_*(\nabla_ST)_{D^\perp}+ F_*\a^j(S,T)
=F_*j_*(\nabla_ST)_{D^\perp}+F_*({}^N\nabla_Sj_*T)_{N_jM}\\
&=F_*({}^N\nabla_Sj_*T)_{j_*D^\perp\oplus N_jM}.
\end{align*}
On the other hand, we have
$$
{}^N\nabla_Sj_*T=({}^N\nabla_Sj_*T)_{j_*D}
+({}^N\nabla_Sj_*T)_{j_*D^\perp\oplus N_jM}.
$$
Since the distribution $\Delta^c$ is totally geodesic and
$j_*D\subset\Delta^c$ it follows that
$$
R=({}^N\nabla_Sj_*T)_{j_*D^\perp\oplus N_jM}\in\Gamma(j^*\Delta^c).
$$
The restriction of $F$ to the leaf of $\Delta^c$ that contains 
$j(\pi(\lambda))$ is an open subset of the affine subspace 
$F(j(\pi(\lambda)))+F_*\Delta^c(j(\pi(\lambda)))$ of $\R^{2n+p}$. 
Then $\bar{F}(\lambda)=F(j(\pi(\lambda)))+F_*R\in F(N)$ as we wished.
\vspace{2ex}\qed

The following result generalizes the Theorem $2$ in \cite{DG4} 
since our rank hypothesis is weaker than the corresponding 
assumption there.

\begin{theorem}\po
Let $f\colon M^{2n}\to\R^{2n+3}$, $n\geq 4$, be a real Kaehler 
submanifold with rank $\varrho_f>3$ everywhere. Then $f$ is locally 
isometrically rigid unless there exists an open subset
$U\subset M^{2n}$ such that the Kaehler extension $F$ of $f|_U$ 
is either a flat or a minimal hypersurface. In this case, any isometric 
deformation of $f|_U$ is the restriction of an isometric 
deformation of $F$.
\end{theorem}

\proof  By  Lemma \ref{unique} and Lemma \ref{c1c2} let 
$V\subset M^{2n}$ be an open subset of $M^{2n}$ on which there 
are a unique pair $(L,\J)$ with $L\subset N_1$ of rank $2$ and 
a tangent subbundle $D=\mathcal{N}_c(\a_{L^\perp})$ with 
$\rank D\geq 2n-2$ such that the conditions $(\mathcal{C}_1)$ 
to $(\mathcal{C}_3)$ are satisfied. 

The subset $U\subset V$ of points where $\dim D=2n-2$ is open
and dense since otherwise $f$ would not be locally substantial.
We argue for points in $U$. Since
$\mathcal{N}_c(\a)=\mathcal{N}_c(\a_L)\cap D$ and 
$\dim\mathcal{N}_c(\a)\leq 2n-8$ by the rank assumption, then
\be\label{complex}
\dim(\mathcal{N}_c(\a_L)+D)=\dim\mathcal{N}_c(\a_L)
+\dim D-\dim\mathcal{N}_c(\a)\geq\dim\mathcal{N}_c(\a_L)+6.
\ee
Thus $\dim\mathcal{N}_c(\a_L)\leq 2n-6$.  Since 
$\mathcal{N}(\a_L)=\mathcal{N}_c(\a_L)$ by \eqref{Jsecondf}
hence $\dim\mathcal{N}(\a_L)\leq 2n-6$.

We claim that $\dim\mathcal{N}(\a)<2n-6$ on $U$.  Since
$\mathcal{N}(\a)\subset\mathcal{N}(\a_L)$ we may 
suppose that $\dim\mathcal{N}(\a_L)=2n-6$. Since the subspace 
$\mathcal{N}_c(\a_L)+D$ is $J$-invariant then either 
$\dim\mathcal{N}(\a_L)+D=2n-2$ or $\mathcal{N}(\a_L)+D=TM$. 
In the former case, the equality part in \eqref{complex} gives
$\varrho_f=3$, a contradiction. Hence $\mathcal{N}(\a_L)+D=TM$. 
In that case and since $\dim D=2n-2$ then
$\dim\mathcal{N}(\a_{L^\perp})\leq 2n-1$.
Hence, from
$\mathcal{N}(\a)=\mathcal{N}(\a_L)\cap\mathcal{N}(\a_{L^\perp})$
we obtain
\begin{equation*}
2n=\dim(\mathcal{N}(\a_L)+\mathcal{N}(\a_{L^\perp}))
=\dim\mathcal{N}(\a_L)+\dim\mathcal{N}(\a_{L^\perp})
-\dim\mathcal{N}(\a)
\end{equation*}
and thus $\dim\mathcal{N}(\a)<2n-6$ as claimed.
Now the proof follows from Theorem $2$ in \cite{DG4}. \qed

\section{Appendix}

The local structures of the substantial real Kaehler submanifolds 
$F\colon M^{2n}\to\R^{2n+2}$ for $n\geq 3$ are discussed next
with separation in the cases were $F$ is a minimal submanifold 
or is free of points where it is minimal. 
\vspace{1ex}
 
\noindent\emph{The non-minimal case}. We assume further that 
either $M^{2n}$ is flat or nowhere flat. In the non flat case the
classifications given below was obtained from \cite{FZ2}.
\vspace{1ex}

\noindent$(i)$ If $M^{2n}=U\subset\R^{2n+2}$ is an open subset
where either $N_1=\Sal(\a^F)$ has rank one  and it is 
not parallel in the normal connection at any point or it
satisfies that rank $N_1=2$ everywhere.
\begin{itemize}
\item[(a)] If rank $N_1=1$ then Theorem $1$ in \cite{DT1} 
gives that $F=H\circ i$ where $i\colon U\to V$ is a totally 
geodesic inclusion with $V\subset\R^{2n+1}$ an open subset and 
$H\colon V\to\R^{2n+2}$ an isometric immersion free of totally 
geodesic points. Moreover, if $Z(i(y))\in T_{i(y)}V$ is an eigenvector 
corresponding to the unique nonzero principal curvature of $H$ then 
the conditions $Z(i(y))\not\in i_*T_yU$ and  $Z(i(y))\not\in N_{i(y)}U$ 
hold at any $y\in U$. 
\item[(b)] If rank $N_1=2$ the nicest local parametric 
classification is given by Corollary $18$ in \cite{FF}.
\end{itemize}

\noindent $(ii)$  An open subset of a cylinder  
$h\times id\colon L^2\times\C^{n-1}\to\R^{2n+2}$ over a non flat 
and nowhere minimal surface $h\colon L^2\to\R^4$. 
\vspace{1ex}

\noindent $(iii)$ A composition of isometric immersion $F=h\circ f$  
where $f\colon M^{2n}\to V\subset\R^{2n+1}$ is a non flat real Kaehler 
hypersurface and $h\colon V\to\R^{2n+2}$ is not totally geodesic.
\vspace{1ex}

\noindent$(iv)$ An open subset of an extrinsic product of two 
Euclidean real Kaehler hypersurfaces where at least one is not 
neither flat nor minimal. 
\vspace{2ex}

\noindent\emph{The minimal case}. For minimal real Kaehler 
submanifolds $F\colon M^{2n}\to\R^N$ in any codimension there 
is the representation given in \cite{APS} and discussed in the 
Appendix of Chapter~$15$ in \cite{DT}.  
Very roughly, the holomorphic representative of the submanifold is 
determined by a set of $n$ independent holomorphic functions  
which span an isotropic subspace of $\C^N$ and have to satisfy 
certain integrability conditions, thus this cannot be seen as a 
classification.
\vspace{1ex}

\noindent $(i)$ If $\varrho_F=1$ we have:
\begin{itemize}
\item[(a)]  An open subset of a cylinder  
$h\times id\colon L^{2}\times\C^{n-}\to\R^{2n+2}$
over a substantial minimal surface $h\colon L^2\to\R^4$.

\item[(b)] A cylinder over a submanifold  parametrically
classified by Theorem $27$ in \cite{DF} by means 
of a Weierstrass type representation given in terms
of $(m-1)$-isotropic surface. These surfaces have been 
completely described in \cite{DF}.
\end{itemize}

\noindent $(ii)$ If $\varrho_F=2$ we have:
\begin{itemize}
\item[(a)] An open subset of the extrinsic product of two 
minimal Euclidean real Kaehler hypersurfaces.
\item[(b)] Examples can be constructed by the use of the 
representation $(8)$ in \cite{DG3} as explained by part $(i)$ of 
the Remark given there. If $M^{2n}$ is complete there is the 
parametric classification provided in \cite{DG3}. A classification 
in the local case remains an open problem unless $n=2$ for which 
there is the classification obtained in \cite{He}. Finally, 
for complete examples for $n=2$ see \cite{DG4}.
\vspace{1ex}
\end{itemize}

\section*{Acknowledgment}

Marcos Dajczer research is part of the project PGC2018-097046-B-I00,
supported by\\ MCIN/AEI/10.13039/501100011033/ FEDER ``Una manera de 
hacer Europa"
\vspace{1ex}

Marcos Dajczer thanks the Mathematics Department of the University
of Murcia where this work was developed for the kind hospitality
during their visit.

{\renewcommand{\baselinestretch}{1}
\hspace*{-30ex}\begin{tabbing}
\indent \= Sergio Julio Chion Aguirre\\
\>CENTRUM Catolica Graduate Business School, Lima, Peru\\
\>Pontificia Universidad Catolica del Peru, Lima, Peru,\\
\>sjchiona@pucp.edu.pe 
\end{tabbing}}

{\renewcommand{\baselinestretch}{1}
\hspace*{-30ex}\begin{tabbing}
\indent \= Marcos Dajczer\\
\>IMPA -- Estrada Dona Castorina, 110\\
\> 22460-320,\\
\>Rio de Janeiro -- Brazil\\
\> marcos@impa.br
\end{tabbing}}

\end{document}